\documentclass[a4paper,9pt]{extarticle}

\usepackage{amsmath,amsthm,amssymb}
\usepackage[authoryear]{natbib}
\usepackage[final]{showkeys}
\usepackage{microtype}
\usepackage{graphicx}
\usepackage{bm}
\usepackage{setspace}

\allowdisplaybreaks
\numberwithin{equation}{section}
\newtheorem{thm}{Theorem}[section]
\newtheorem{rem}{Remark}[section]

\usepackage[hmargin=2.5cm,vmargin={3cm,3.5cm}]{geometry}

\title{On the Integral of Fractional Poisson Processes}

\author{$\text{Enzo Orsingher}_1$, $\text{Federico Polito}_2$\\
	\footnotesize (1) -- Dipartimento di Scienze Statistiche, ``Sapienza'' Universit\`a di Roma\\
	\footnotesize Email address: enzo.orsingher@uniroma1.it\\
	\footnotesize (2) -- Dipartimento di Matematica, Universit\`a di Torino\\
	\footnotesize Tel: +39 011 6702937, fax: +39 011 6702878\\
	\footnotesize Email address: federico.polito@unito.it (Corresponding author)}

\begin{document}

\onehalfspacing

\maketitle

\begin{abstract}

	In this paper we consider the Riemann--Liouville fractional integral $\mathcal{N}^{\alpha,\nu}(t)=
	\frac{1}{\Gamma(\alpha)} \int_0^t (t-s)^{\alpha-1}N^\nu(s) \, \mathrm ds $, where
	$N^\nu(t)$, $t \ge 0$,
	is a fractional Poisson process of order $\nu \in (0,1]$, and $\alpha > 0$.
	We give the explicit bivariate distribution
	$\Pr \{ N^\nu(s)=k, N^\nu(t)=r \}$, for $t \ge s$, $r \ge k$, the mean
	$\mathbb{E}\, \mathcal{N}^{\alpha,\nu}(t)$ and the variance $\mathbb{V}\text{ar}\,  \mathcal{N}^{\alpha,\nu}(t)$.
	We study the process $\mathcal{N}^{\alpha,1}(t)$ for which we are able to produce explicit results
	for the conditional and absolute variances and means.
	Much more involved results on $\mathcal{N}^{1,1}(t)$ are presented in the last section where also
	distributional properties of the integrated Poisson process
	(including the representation as random sums) is derived.
	The integral of powers of the Poisson process is examined and its connections with generalised
	harmonic numbers is discussed.
	
	\vspace{.2cm}
	

\end{abstract}

	\section{Introduction}
	
		The fractional Poisson process $N^\nu(t)$, $t \ge 0$, $0<\nu \le 1$, has been introduced and studied in the
		last decade by \citet{laskin}, \citet{mainardi}, \citet{beg}, \citet{scalas}.
		The starting point of the investigations of some authors was the derivation of the distribution
		\begin{align}
			\label{00process}
			p_k^\nu(t) = \Pr \{ N^\nu(t) = k \}, \qquad k \ge 0,
		\end{align}
		by solving the fractional equations
		\begin{align}
			\label{ccca}
			\begin{cases}
				\frac{\mathrm d^\nu}{\mathrm dt^\nu} p_k^\nu(t) = -\lambda p_k^\nu(t) + \lambda p_{k-1}^\nu(t), \\
				p_k^\nu(0) =
				\begin{cases}
					1, & k = 0, \\
					0, & k \ge 1.
				\end{cases}
			\end{cases}
		\end{align}

		The derivative appearing in \eqref{ccca} is meant in the sense of Riemann--Liouville in
		\citet{laskin} and in the sense of Dzhrbashyan--Caputo in \citet{beg}. The distribution \eqref{00process}
		reads
		\begin{align}
			\label{00tex}
			p_k^\nu(t) & = \sum_{r=k}^\infty (-1)^{r-k} \binom{r}{k} \frac{(\lambda t^\nu)^r}{\Gamma(\nu r+1)}
			=\frac{(\lambda t^\nu)^k}{k!} \sum_{r=0}^\infty \frac{(r+k)!}{r!}
			\frac{(-\lambda t^\nu)^r}{\Gamma(\nu(k+r)+1)}.
		\end{align}

		The fractional Poisson process is also constructed as a renewal process in \citet{mainardi}
		and \citet{beg} where is shown that its distribution coincides with \eqref{00tex}.
		\citet{meers} treat and analyze in a unified way the process obtained from the time-fractional
		equation and the related renewal process with Mittag--Leffler distributed interarrival times.
		Furthermore, other generalizations in a fractional sense or as a renewal process
		with generalized Mittag--Leffler waiting times have been introduced by \citet{space}
		and \citet{caho}.
		
		The intertime $T_1^\nu$ between successive events has distribution
		\begin{align}
			\label{nnn}
			\Pr \{ T_1^\nu \in \mathrm ds \} = \lambda s^{\nu-1} E_{\nu,\nu} (-\lambda s^\nu) \, \mathrm ds,
			\qquad s \ge 0,
		\end{align}
		while the waiting time for the $k$-th event $T_k^\nu$ has distribution
		\begin{align}
			\label{00cc}
			\Pr \{ T_k^\nu \in \mathrm ds \} = \lambda^k s^{\nu k-1} E_{\nu,\nu k}^k (-\lambda s^\nu)
			\, \mathrm ds, \qquad s \ge 0,
		\end{align}
		where
		\begin{align}
			E_{\alpha,\eta}^\gamma (z) = \sum_{r=0}^\infty \frac{(\gamma)_r z^r}{r! \Gamma(\alpha r+\eta)},
			\quad \alpha,\eta,\gamma \in \mathbb{C}, \: \Re(\alpha),\Re(\eta),\Re(\gamma)>0, \: z \in \mathbb{C},
		\end{align}
		is the generalised Mittag--Leffler function \citep{sax}.
		Note that $(\gamma)_r = \gamma(\gamma+1)\dots (\gamma+r-1)$,
		$\gamma \ne 0$, and that $E_{\alpha,\eta}(z) = E_{\alpha,\eta}^1(z)$.
		
		The multivariate distribution of the fractional Poisson process
		$
			\Pr \{ N^\nu(t_1) = n_1, \dots, N^\nu(t_k) = n_k \},
		$
		where $t_1<t_2<\dots <t_k$, $n_1\le n_2 \le \dots \le n_k$, can be written down by considering its renewal
		structure and by exploiting \eqref{00cc} and \eqref{nnn} (see e.g.\ \citet{scalas}).
		We are able to give the explicite bivariate distribution
		in terms of generalised
		Mittag--Leffler functions. This plays a crucial role in the analysis of the variance
		of the fractional integral of the fractional Poisson process, i.e.\
		\begin{align}
			\label{bababa}
			\mathcal{N}^{\alpha,\nu}(t) = \frac{1}{\Gamma(\alpha)} \int_0^t (t-s)^{\alpha-1} N^\nu(s) \, \mathrm ds,
			\qquad t \ge 0, \: 0<\nu\le 1, \: \alpha > 0,
		\end{align}
		with
$
			\mathbb{E} \, \mathcal{N}^{\alpha,\nu} (t) = \lambda t^{\alpha+\nu} / \Gamma(\alpha+\nu+1).
$
		
		For $\nu=1$, we obtain the fractional integral of the classical Poisson process
		with intensity $\lambda$.
		The motivation in studying the above process is based on the fact that
		integrated non-negative stochastic processes and in general	integrated counting processes
		often arise in the applied mathematical literature (see for example \citet{jerwood},
		\citet{downton}, \citet{hernandez}, \citet{stefanov}, \citet{pollett}, and
		the references therein). The analysis of the integrated process \eqref{bababa} is interesting
		as it permits to generalize the ideas behind such studies to a non-integer framework.
		Note also that for $\alpha \in \mathbb{N}$
		the Riemann--Liouville fractional integral coincides with a classical multiple integral.
		
		The main result
		for the Riemann--Liouville integral is the conditional second moment
		\begin{align}
			\label{00lab}
			\mathbb{E} \left\{ \left( \frac{1}{\Gamma(\alpha)} \int_0^t (t-s)^{\alpha-1} N(s) \, \mathrm ds
			\right)^2 \biggr| N(t) = n \right\}
			= \frac{2nt^{2\alpha}\Gamma(2\alpha)}{\alpha\Gamma^2(\alpha)\Gamma(2\alpha+2)}
			+\frac{n(n-1)t^{2\alpha}}{\Gamma^2(\alpha+2)}.
		\end{align}
		Of course we have that \citep[page 21]{kingman}
		\begin{align}
			\label{00lab2}
			\mathbb{E} \left\{ \frac{1}{\Gamma(\alpha)} \int_0^t (t-s)^{\alpha-1} N(s) \, \mathrm ds
			\biggr| N(t) = n \right\} = \frac{n t^\alpha}{\Gamma(\alpha+2)},
		\end{align}
		and thus
		\begin{align}
			\mathbb{E} \left\{ \frac{1}{\Gamma(\alpha)}
			\int_0^t (t-s)^{\alpha-1} N(s) \, \mathrm ds \right\} =
			\frac{\lambda t^{\alpha+1}}{\Gamma(\alpha+2)}.
		\end{align}
		In light of \eqref{00lab} and \eqref{00lab2} we are able to give the conditional variance
		of the fractional integral of the Poisson process:
		\begin{align}
			\label{00ba}
			\mathbb{V}\text{ar} \left\{ \left. \frac{1}{\Gamma(\alpha)} \int_0^t (t-s)^{\alpha-1} N(s) \, \mathrm ds
			\right| N(t) = n \right\}
			= \frac{n t^{2\alpha}\alpha^2}{(2\alpha+1)\Gamma^2(\alpha+2)}.
		\end{align}
		Therefore we extract from \eqref{00ba} and \eqref{00lab2} that
		\begin{align}
			\label{06finintro}
			\mathbb{V}\text{ar} \left( \mathcal{N}^{\alpha,1}(t) \right) =
			\frac{\lambda t^{2\alpha+1}}{(2\alpha+1)\Gamma^2(\alpha+1)}.
		\end{align}
		
		For $\alpha=1$ we have the integral of the classical Poisson process which can be written as a
		random sum, i.e.\
		\begin{align}
			\label{00fr}
			\int_0^t N(s) \, \mathrm ds \overset{\text{d}}{=} \sum_{j=1}^{N(t)} X_j.
		\end{align}
		The random variables $X_j$s appearing in \eqref{00fr} are i.i.d.\ with uniform law in $(0,t)$
		independent of $N(t)$.
		In \eqref{00fr} we consider that the sum on the right hand side is empty for $N(t)=0$.
		For the conditional integral of the Poisson process we have that
		\begin{align}
			\label{00qa}
			\mathbb{E} \left\{ \int_0^t N(s) \, \mathrm ds \biggr| N(t) = n \right\}
			= \frac{nt}{2}, \qquad
			\mathbb{V}\text{ar} \left\{ \int_0^t N(s) \, \mathrm ds \biggr| N(t) = n \right\}
			= \frac{nt^2}{12},
		\end{align}
		which are also special cases of \eqref{00lab2} and \eqref{00ba} for $\alpha=1$. The results
		\eqref{00qa} are also obtained by a different, alternative method.
		
		Finally we examine in Section \ref{prop99} the process
		$
			\mathring{N} (t) = N_\lambda(t) - N_\beta(t)
		$, $t \ge 0$,
		where $N_\lambda(t)$ and $N_\beta(t)$ are independent Poisson processes of parameter
		$\lambda>0$ and $\beta>0$, respectively.
		It is well-known that
		\begin{align}
			\label{00gin}
			\Pr \{ \mathring{N}(t) = r  \} = e^{-(\beta+\lambda)t} \left( \frac{\lambda}{\beta} \right)^{r/2}
			I_{|r|}(2 t \sqrt{\lambda \beta}), \qquad r \in \mathbb{Z}, \: t\ge 0,
		\end{align}
		where
		\begin{align}
			I_\xi(z) = \sum_{k=0}^\infty \left( \frac{z}{2} \right)^{2k+\xi} \frac{1}{k!\Gamma(k+\xi+1)}
		\end{align}
		is the modified Bessel function of the first kind.
		The distribution \eqref{00gin} is called the Skellam distribution.
		For the integral process
		\begin{align}
			\int_0^t \mathring{N}(s) \, \mathrm ds = \int_0^t N_\lambda(s) \, \mathrm ds
			- \int_0^t N_\beta(s) \, \mathrm ds,
		\end{align}
		we show that
		\begin{align}
			\label{08ciao}
			\int_0^t \mathring{N}(s) \, \mathrm ds \overset{\text{d}}{=} \sum_{j=1}^{\widetilde{N}(t)} Z_j,
		\end{align}
		where $\widetilde{N}(t)$, $t \ge 0$, is a Poisson process with rate $\lambda+\beta$ and
		the random variables $Z_j$s are i.i.d.\ with density
		\begin{align}
			\label{00io}
			f(s) =
			\begin{cases}
				\frac{\beta}{t(\lambda+\beta)}, & -t < s \le 0, \\
				\frac{\lambda}{t(\lambda+\beta)}, & 0 < s < t.
			\end{cases}
		\end{align}
		Clearly, for $\beta=\lambda$, \eqref{00io} becomes the uniform distribution in $(-t,t)$.
		As before, in \eqref{08ciao} the sum is considered equal to zero when $\widetilde{N}(t) = 0$.
		We remark that integrals of different point processes have been considered over the years,
		for example in \citet{puri},
		where the integral of the birth and death process has been examined.

	\section{Fractional integral of the fractional Poisson process}
	
		For the fractional Poisson process $N^\nu(t)$, $t \ge 0$, described in the introduction
		we consider the Riemann--Liouville fractional integral
		\begin{align}
			\label{00car}
			\mathcal{N}^{\alpha,\nu}(t) = \frac{1}{\Gamma(\alpha)} \int_0^t (t-s)^{\alpha-1} N^\nu(s) \, \mathrm ds,
			\qquad t \ge 0, \: 0<\nu\le 1, \: \alpha > 0.
		\end{align}
		For integer values of $\alpha$, say $\alpha = m$, the integral \eqref{00car} can be written as
		\begin{align}
			& \frac{1}{(m-1)!} \int_0^t (t-s)^{m-1} N^\nu(s) \, \mathrm ds
			= \int_0^t \mathrm ds_1 \int_{s_1}^t \mathrm ds_2 \dots \int_{s_{m-1}}^t 
			N^\nu(s_m)\, \mathrm ds_{m} ,
		\end{align}
		and this offers an intuitive interpretation of \eqref{00car}. By taking into account
		formula (2.7) of \citet{beg} it is a quick matter to check that
		\begin{align}
			\label{02comp}
			\mathbb{E} \, \mathcal{N}^{\alpha,\nu} (t) & =\frac{1}{\Gamma(\alpha)} \int_0^t (t-s)^{\alpha-1}
			\mathbb{E}N^\nu(s) \, \mathrm ds
			= \frac{1}{\Gamma(\alpha)} \int_0^t (t-s)^{\alpha-1} \frac{\lambda s^\nu}{\Gamma(\nu+1)} \, \mathrm ds
			= \frac{\lambda t^{\alpha+\nu}}{\Gamma(\alpha+\nu+1)}.
		\end{align}
		Note that if $0 < \alpha+\nu \le 1$, then $\mathbb{E} \, \mathcal{N}^{\alpha,\nu} (t)
		= \mathbb{E} N^{\alpha+\nu}(t)$.
		
		The fractional Poisson process can be seen as a renewal process with intertime between successive
		events possessing distribution
		\begin{align}
			\label{fufu}
			\Pr \{ T_1^\nu \in \mathrm ds \} = \lambda s^{\nu-1} E_{\nu,\nu}(-\lambda s^\nu) \, \mathrm ds,
			\qquad s \ge 0, \: 0<\nu \le 1.
		\end{align}
		This has been proved in \citet{mainardi}, \citet{beg}, and \citet{scalas}.
		The random instant of the occurrence of the $k$th event for $N^\nu(t)$, $t \ge 0$, is denoted by
		$
			T_k^\nu = \inf \{ t \colon N^\nu(t) = k \}.
		$
		We need also the symbol
		$
			\mathcal{T}_h^{\nu,k} = T^\nu_{k+h} - T^\nu_k \overset{\text{d}}{=} T_h^\nu,
		$
		where $\mathcal{T}_h^{\nu,k}$ represents the length of the time interval separating the
		$k$th and the $(k+h)$th event. The distribution of $T_k^\nu$ is given in \citet{beg2}
		as
		\begin{align}
			\label{fufu1}
			\Pr \{ T_k^\nu \in \mathrm ds \} = \lambda^k s^{\nu k-1} E_{\nu,\nu k}^k (-\lambda s^\nu)\, \mathrm ds,
			\qquad s \ge 0, \: 0 < \nu \le 1.
		\end{align}
		
		\begin{thm}
			The bivariate distribution of the fractional Poisson process reads
			\begin{align}
				\label{01dist}
				\Pr & \{ N^\nu(s) = k, N^\nu(t) = r\} \\
				= {} & \lambda^r \int_0^s w^{\nu k-1} E_{\nu,\nu k}^k (-\lambda w^\nu) \, \mathrm dw \int_{s-w}^{t-w}
				y^{\nu-1} E_{\nu,\nu}(-\lambda y^\nu) (t-w-y)^{\nu(r-k-1)}
				E_{\nu,\nu(r-k-1)+1}^{r-k}
				(-\lambda(t-w-y)^\nu) \, \mathrm dy, \notag
			\end{align}
			where
			\begin{align}
				E^\gamma_{\alpha,\eta} (z) = \sum_{r=0}^\infty \frac{\Gamma(\gamma +r)z^r}{r!\Gamma(\gamma)
				\Gamma(\alpha r + \eta)},
				\quad \alpha,\eta,\gamma \in \mathbb{C}, \: \Re(\alpha),\Re(\eta),\Re(\gamma)>0, \: z \in \mathbb{C},
			\end{align}
			is the generalised Mittag--Leffler function.
			\begin{proof}
				In order to obtain the distribution \eqref{01dist} we must have a look at Figure \ref{01fig}, where
				the instants of occurrence of the relevant events are depicted.
				\begin{figure}
					\centering
					\includegraphics[scale=1]{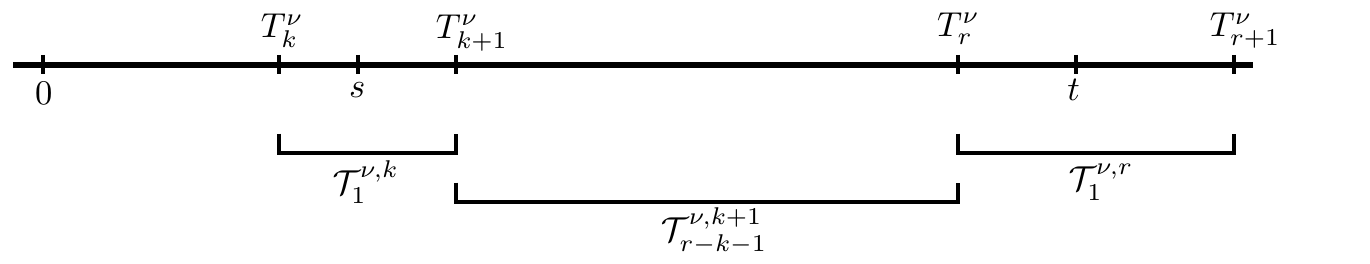}
					\caption{\label{01fig}The instants of occurrence of the events and the related
						waiting times.}
				\end{figure}
				The bivariate distribution can be written as
				\begin{align}
					\Pr & \{ N^\nu(s) = k, N^\nu(t) = r\}
					= \iiiint_D \Pr \{ T_k^\nu \in \mathrm dw, \mathcal{T}_1^{\nu,k} \in \mathrm dy,
					\mathcal{T}_{r-k-1}^{\nu,k+1} \in \mathrm d\xi, \mathcal{T}_1^{\nu,r} \in \mathrm d\eta \},
				\end{align}
				where
				$
					D = \{ (0<w<s) \cap (y+w>s)
					\cap (t>y+w+\xi>s) \cap (y+w+\xi +\eta) > t \}.
				$
				We have that (by keeping in mind the independence of the intertimes between events)
				\begin{align}
					\label{01bic}
					\Pr & \{ N^\nu(s) = k, N^\nu(t) = r\} \\
					= {} & \int\limits_0^s \int\limits_{s-w}^{t-w} \int\limits_0^{t-(w+y)} \int\limits_{t-(y+w+\xi)}^\infty
					\Pr \{ T_k^\nu \in \mathrm dw, \mathcal{T}_1^{\nu,k} \in \mathrm dy,
					\mathcal{T}_{r-k-1}^{\nu,k+1} \in \mathrm d\xi, \mathcal{T}_1^{\nu,r} \in \mathrm d\eta \} \notag \\
					= {} & \int_0^s \Pr \{ T^\nu_k \in \mathrm dw \} \int_{s-w}^{t-w} \Pr \{ \mathcal{T}^{\nu,k}_1
					\in \mathrm dy \}
					\int_0^{t-(w+y)} \Pr \{ \mathcal{T}_{r-k-1}^{\nu,k+1} \in \mathrm d \xi \}
					\int_{t-(y+w+\xi)}^\infty \Pr \{ \mathcal{T}_1^{\nu,r} \in \mathrm d\eta \} \notag \\
					= {} & \int_0^s \lambda^k w^{\nu k -1} E_{\nu,\nu k}^k (-\lambda w^\nu) \int_{s-w}^{t-w}
					\lambda y^{\nu-1} E_{\nu,\nu}(-\lambda y^\nu)
					\int_0^{t-(w+y)} \lambda^{r-k-1} \xi^{\nu(r-k-1)-1}
					E_{\nu,\nu(r-k-1)}^{r-k-1} (-\lambda \xi^\nu) \notag \\
					& \times E_{\nu,1}(-\lambda (t-y-w-\xi)^\nu) \,
					\mathrm d\xi \, \mathrm dy \, \mathrm dw. \notag
				\end{align}
				Writing down the integral in \eqref{01bic} consider Figure \ref{01fig}
				and the independence of the intertimes $T_k^\nu$, $\mathcal{T}_1^{\nu,k}$,
				$\mathcal{T}_{r-k-1}^{\nu,k+1}$ and $\mathcal{T}_1^{\nu,r}$ (with distributions
				\eqref{fufu} and \eqref{fufu1}).
				
				Formula \eqref{01bic}
				can be further simplified by recurring to the following relation
				(see e.g.\ \citet{sax}, formula (11.7), page 17):
				\begin{align}
					\label{01for}
					\int_0^x (x-t)^{\beta-1} E_{\alpha,\beta}^\gamma [a(x-t)^\alpha] t^{\zeta-1}
					E_{\alpha,\zeta}^\sigma (at^\alpha)
					\mathrm dt = x^{\beta + \zeta-1} E_{\alpha, \beta+\zeta}^{\gamma+\sigma}
					(ax^\alpha),
				\end{align}
				where $\alpha,\beta,\gamma,a,\zeta,\sigma \in \mathbb{C}$, and $\Re(\alpha)>0$, $\Re(\beta)>0$,
				$\Re(\gamma)>0$, $\Re(\zeta)>0$, $\Re(\sigma)>0$.
				With the necessary substitutions in \eqref{01for}, that is,
				$x=t-(w+y)$, $t=\xi$, $\alpha=\nu$, $\zeta=\nu(r-k-1)$, $\sigma= r-k-1$, $a=-\lambda$, and
				$\beta=1$,				
				we have that
				\begin{align}
					\Pr & \{ N^\nu(s) = k, N^\nu(t) = r\} \\
					= {} & \lambda^r \int_0^s w^{\nu k-1} E_{\nu,\nu k}^k (-\lambda w^\nu) \int_{s-w}^{t-w}
					y^{\nu-1} E_{\nu,\nu}(-\lambda y^\nu) (t-w-y)^{\nu(r-k-1)}
					E_{\nu,\nu(r-k-1)+1}^{r-k}
					(-\lambda(t-w-y)^\nu) \, \mathrm dy \, \mathrm dw. \notag
				\end{align}
			\end{proof}
		\end{thm}
		
		\begin{rem}
			We show now that \eqref{01dist}, for $\nu=1$, that is for the classical homogenous Poisson process, yields
			\begin{align}
				\Pr \{ N^1(s)=k,N^1(t) = r \} = \frac{\lambda^r s^k(t-s)^{r-k}}{r!} \binom{r}{k}
				e^{-\lambda t}.
			\end{align}
			Since
			$
				E_{1,1}(x) = e^x
			$ and
			$
				E_{1,k}^k(x) = e^x /(k-1)!
			$
			we have that
			\begin{align}
				\Pr & \{ N^1(s)= k, N^1(t)=r \}
				= \lambda^r \frac{s^k}{k!} \frac{(t-s)^{r-k}}{(r-k)!} e^{-\lambda t}. 
			\end{align}
		\end{rem}
		
		\begin{rem}
			If we change the variable in the outer integral of \eqref{01dist} we get
			\begin{align}
				\label{01new}
				\Pr & \{ N^\nu(s) = k, N^\nu(t) = r\} \\
				= {} & \lambda^r \int_{t-s}^t (t-z)^{\nu k-1} E_{\nu,\nu k}^k (-\lambda (t-z)^\nu) \, \mathrm dz
				\int_{z+s-t}^z
				y^{\nu-1} E_{\nu,\nu}(-\lambda y^\nu) (z-y)^{\nu(r-k-1)}
				E_{\nu,\nu(r-k-1)+1}^{r-k}
				(-\lambda(z-y)^\nu) \, \mathrm dy. \notag
			\end{align}
			In \eqref{01new} we have an integral of the form
			\begin{align}
				\label{01prab}			
				\left( \bm{\mathrm E}^\gamma_{\rho,\mu,\omega;a+} \phi \right)(x) = \int_a^x (x-t)^{\mu-1}
				E_{\rho,\mu}^\gamma
				\left(\omega(x-t)^\rho\right) \phi(t) \, \mathrm dt, \qquad x>a,
			\end{align}
			treated in \citet{prab,saigo,kil,srivastava}.
			The integral in \eqref{01prab} is a generalisation of the Riemann--Liouville fractional integral.
			
			In light of \eqref{01prab}, the bivariate distribution \eqref{01new} can be written as
			\begin{align}
				\Pr & \{ N^\nu(s) = k, N^\nu(t) = r\} \\
				= {} & \lambda^r \int_{t-s}^t (t-z)^{\nu k-1} E_{\nu,\nu k}^k \left( -\lambda(t-z)^\nu \right)
				\left( \bm{\mathrm E}^{r-k}_{\nu,\nu(r-k-1)+1,-\lambda;(z+s-t)+}
				y^{\nu-1}E_{\nu,\nu}(-\lambda y^\nu) \right)
				(z) \, \mathrm dz \notag \\
				= {} & \lambda^r \left( \bm{\mathrm E}^k_{\nu,\nu k,-\lambda;(t-s)+}
				\left( \bm{\mathrm E}^{r-k}_{\nu,\nu(r-k-1)+1,-\lambda;(z+s-t)+}
				 y^{\nu-1}E_{\nu,\nu}(-\lambda y^\nu) \right)
				(z) \right) (t). \notag
			\end{align}
		\end{rem}
		
		\begin{rem}
			For the second-order moment of the fractional integral we have that
			\begin{align}
				\label{01ref}
				\mathbb{E} \left\{ \frac{1}{\Gamma(\alpha)} \int_0^t (t-u)^{\alpha-1} N^\nu(u)
				\, \mathrm du \right\}^2
				= \frac{1}{\Gamma^2(\alpha)} \int_0^t \int_0^t (t-u)^{\alpha-1}(t-v)^{\alpha-1}
				\mathbb{E} \{ N^\nu(u) N^\nu(v) \} \, \mathrm du \, \mathrm dw,
			\end{align}
			where
			$
				\mathbb{E} \{ N^\nu(u) N^\nu(v) \} = \sum_{k=0}^\infty \sum_{r=k}^\infty
				k \, r \Pr \{ N^\nu(u) = k, N^\nu(v) = r \}
			$, $v> u$.
			Unfortunately the complicated structure of the distribution \eqref{01dist} does not permit us
			to determine a closed form for \eqref{01ref}.
		\end{rem}

	\section{Fractional integral of the homogenous Poisson process}
	
		We now restrict the analysis of \eqref{00car} to the case $\nu=1$, that is we study the Riemann--Liouville
		fractional integral of the homogenous Poisson process:
		\begin{align}
			\label{02def}
			\mathcal{N}^{\alpha,1} (t) = \mathcal{N}^\alpha (t) = \frac{1}{\Gamma(\alpha)}
			\int_0^t (t-s)^{\alpha-1} N(s) \, \mathrm ds,
		\end{align}
		for $\alpha>0$, $t \ge 0$.
		From \eqref{02comp} we have that
		\begin{align}
			\label{02mean}
			\mathbb{E} \, \mathcal{N}^\alpha(t) = \frac{\lambda t^{\alpha+1}}{\Gamma(\alpha+2)}.
		\end{align}
		We are now in the position of evaluating explicitly $\mathbb{E}\left[ \mathcal{N}^\alpha(t) \right]^2$ 
		and also
		$
			\mathbb{E} \left\{ \left[ \mathcal{N}^\alpha(t) \right]^2 | N(t) = n \right\} .
		$
		Thus we can state the following theorem.
		
		\begin{thm}
			The variance of \eqref{02def} has the following form:
			\begin{align}
				\label{05cc}
				\mathbb{V}\text{ar} \, \mathcal{N}^\alpha(t) = \frac{\lambda t^{2\alpha+1}}{(2\alpha+1)\Gamma^2(\alpha+1)},
				\qquad t \ge 0, \: \alpha > 0.
			\end{align} 
			\begin{proof}
				We start by evaluating $\mathbb{E}\left[ \mathcal{N}^\alpha(t) \right]^2$.
				\begin{align}
					\mathbb{E}\left[ \mathcal{N}^\alpha(t) \right]^2 = {} &
					\frac{1}{\Gamma^2(\alpha)} \int_0^t \int_0^t (t-s)^{\alpha-1}(t-w)^{\alpha-1}
					\mathbb{E}\left[ N(s)N(w) \right] \, \mathrm ds \, \mathrm dw \\
					= {} & \frac{2}{\Gamma^2(\alpha)} \int_0^t \int_s^t (t-s)^{\alpha-1} (t-w)^{\alpha-1}
					(\lambda s+\lambda^2 sw) \, \mathrm dw \, \mathrm ds \notag \\
					= {} & \frac{2}{\Gamma^2(\alpha)} \left\{ \lambda \int_0^t s(t-s)^{\alpha-1}
					\mathrm ds \int_s^t (t-w)^{\alpha-1} \mathrm dw \right.
					+ \left. \lambda^2 \int_0^t s(t-s)^{\alpha-1} \mathrm ds \int_s^t w(t-w)^{\alpha-1}
					\mathrm dw \right\} \notag \\
					= {} & \frac{2}{\Gamma^2(\alpha)} \left\{ \frac{\lambda}{\alpha}
					\int_0^t s(t-s)^{2\alpha-1} \mathrm ds + \frac{\lambda^2}{\alpha}
					\int_0^t s^2(t-s)^{2\alpha-1} \mathrm ds \right.
					+ \left. \frac{\lambda^2}{\alpha(\alpha+1)} \int_0^t s(t-s)^{2\alpha} \mathrm ds \right\} \notag \\
					= {} & \frac{2}{\Gamma^2(\alpha)} \left\{ \frac{\lambda t^{2\alpha+1}}{2\alpha^2 (2\alpha+1)}
					+ \frac{\lambda^2 t^{2\alpha +2}}{2\alpha(\alpha+1)^2(2\alpha+1)} \right.
					+ \left. \frac{\lambda^2 t^{2\alpha+2}}{2\alpha^2 (2\alpha+1)(\alpha+1)} \right\} \notag \\
					= {} & \frac{1}{\Gamma^2(\alpha)} \left\{
					\frac{\lambda t^{2\alpha+1}}{\alpha^2(2\alpha+1)} + \frac{\lambda^2
					t^{2\alpha+2}}{\alpha^2(\alpha+1)^2} \right\}
					= \frac{\lambda t^{2\alpha+1}}{(2\alpha+1)\Gamma^2(\alpha+1)} + \frac{\lambda^2
					t^{2\alpha+2}}{\Gamma^2(\alpha+2)}. \notag
				\end{align}
				By considering \eqref{02mean} we immediately arrive at the claimed result.
			\end{proof}
		\end{thm}
		
		\begin{rem}
			For the conditional mean, we directly arrive at the result
			\begin{align}
				\label{02blabla}
				\mathbb{E} \left( \mathcal{N}^\alpha(t)
				\bigr| N(t) = n \right) = \frac{1}{\Gamma(\alpha)} \int_0^t (t-s)^{\alpha-1} \frac{ns}{t} \mathrm ds
				= \frac{n t^\alpha}{\Gamma(\alpha+2)}.
			\end{align}
			In \eqref{02blabla} we considered that
			\begin{align}
				\label{labello}
				\Pr \{ N(s)= r | N(t) = n \} = \binom{n}{r} \left( \frac{s}{t} \right)^r \left( 1-\frac{s}{t}
				\right)^{n-r}, \qquad 0\le r \le n, \: s<t.
			\end{align}			
		\end{rem}
		
		In order to obtain the conditional variance of $\mathcal{N}^\alpha(t)$ we need the following result.
		
		\begin{thm}
			For the homogenous Poisson process we have that
			\begin{align}
				\label{03cut}
				\mathbb{E} \left\{ N(s)N(w) \bigr| N(t)=n \right\} =
				\frac{ns}{t} +n(n-1) \frac{sw}{t^2}, \qquad 0<s<w<t.
			\end{align}
			\begin{proof}
				In order to obtain \eqref{03cut} we evaluate the following bivariate conditional distribution.
				For $s<w<t$ we have that
				\begin{align}
					\label{03cut2}
					\Pr \{ N(s) = h, N(w) = k | N(t) = n \}
					= \Pr \{ N(s) = h | N(w) = k \} \Pr \{ N(w)= k| N(t) = n \}
				\end{align}
				because the time-reversed Poisson process is Markovian.
				From \eqref{03cut} we obtain the following trinomial distribution.
				\begin{align}
					\label{03cut3}
					& \Pr \{ N(s) = h, N(w) = k | N(t) = n \} \\
					& = \binom{k}{h} \left( \frac{s}{w} \right)^h \left( 1-\frac{s}{w} \right)^{k-h}
					\binom{n}{k} \left( \frac{w}{t} \right)^k \left( 1-\frac{w}{t} \right)^{n-k}
					= \frac{n!s^h(w-s)^{k-h}(t-w)^{n-k}}{h!(k-h)!(n-k)! t^n}, \qquad
					0<s<w<t. \notag
				\end{align}
				We evaluate directly the conditional mixed moment of \eqref{03cut3} as follows.
				\begin{align}
					\label{03cc}
					\mathbb{E} & \left\{ N(s)N(w) | N(t) = n \right\}
					= \sum_{h=0}^n \sum_{k=h}^n k \, h \frac{n!s^h(w-s)^{k-h}(t-w)^{n-k}}{h!(k-h)!(n-k)! t^n} \\
					= {} & \frac{n!}{t^n} \sum_{h=1}^n \frac{s^h}{(h-1)!} \sum_{r=0}^{n-h} \frac{(h+r)}{r!
					(n-r-h)!} (w-s)^r (t-w)^{n-r-h} \notag \\
					= {} & \frac{n!}{t^n} \sum_{h=1}^n \frac{s^h}{(h-1)!} \left[
					\frac{h}{(n-h)!} \sum_{r=0}^{n-h} \binom{n-h}{r} (w-s)^r (t-w)^{n-r-h} \right. \notag \\
					& + \left. \frac{(w-s)}{(n-h-1)!} \sum_{l=0}^{n-h-1} \binom{n-h-1}{l}
					(w-s)^l (t-w)^{n-h-1-l} \right] \notag \\
					= {} & \frac{n!}{t^n} \sum_{h=1}^n \frac{s^h}{(h-1)!}
					\left[ \frac{h}{(n-h)!} (t-s)^{n-h} + \frac{(w-s)}{(n-h-1)!}
					(t-s)^{n-1-h} \right] \notag \\
					= {} & \frac{n!}{t^n} \left[ \sum_{m=0}^{n-1} \binom{n-1}{m} s^{m+1}
					(t-s)^{n-1-m} \frac{m+1}{(n-1)!} \right.
					+ \left. \sum_{h=1}^{n-1} \frac{s^h}{(h-1)!} \frac{(w-s)}{(n-h-1)!} (t-s)^{n-1-h} \right] \notag \\
					= {} & \frac{n!}{t^n} \left[ \frac{1}{(n-1)!}
					\left( t^{n-1} (n-1)\frac{s^2}{t} + t^{n-1}s \right)
					+ \frac{(w-s)st^{n-2}}{(n-2)!} \right] \notag \\
					= {} & n(n-1) \frac{s^2}{t^2} + \frac{ns}{t} + (w-s) \frac{s(n-1)n}{t^2}
					= \frac{n(n-1)}{t^2} ws + \frac{ns}{t}. \notag
				\end{align}
			\end{proof}
		\end{thm}
		
		\begin{rem}
			As a simple check we note that
			\begin{align}
				\mathbb{E} \left[ \mathbb{E} \left( N(s)N(w)|N(t) \right) \right] =
				\lambda^2 ws + \lambda s, \qquad 0<s<t.
			\end{align}
			Since for $w=s$
			\begin{align}
				\mathbb{E} \left( [N(s)]^2 | N(t) = n \right) = \frac{n(n-1)}{t^2}s^2 + \frac{ns}{t},
			\end{align}
			the conditional variance reads
			\begin{align}
				\mathbb{V}\text{ar} \left( N(s)|N(t) = n \right) = \frac{ns}{t} -\frac{ns^2}{t^2}.
			\end{align}
			In turn, the unconditional variance can be obtained as follows.
			\begin{align}
				\mathbb{V}\text{ar} \, N(s) & = \mathbb{E} \left[ \mathbb{V}\text{ar} \left( N(s)|N(t) \right) \right]
				+ \mathbb{V}\text{ar} \left[ \mathbb{E} \left( N(s)|N(t) \right) \right] \\
				& = \left( \frac{s}{t}-\frac{s^2}{t^2} \right) \mathbb{E} N(t) + \mathbb{V}\text{ar}
				\left( \frac{s}{t}N(t) \right)
				= \left( \frac{s}{t}-\frac{s^2}{t^2} \right) \lambda t + \frac{s^2}{t}\lambda
				= \lambda s \notag.
			\end{align}
		\end{rem}
		
		We arrive at the conditional variance of the Riemann--Liouville fractional integral of the Poisson
		process in the next theorem.
		
		\begin{thm}
			\label{06teorema}
			We have that
			\begin{align}
				\label{05dd}
				\mathbb{V}\text{ar} \left( \mathcal{N}^\alpha(t) | N(t) = n \right)
				= \frac{nt^{2\alpha}\alpha^2}{(2\alpha+1)\Gamma^2(\alpha+2)}.
			\end{align}
			\begin{proof}
				Exploiting result \eqref{03cc} we have
				\begin{align}
					\label{03re}
					\mathbb{E} & \left\{ \left. \left( \frac{1}{\Gamma(\alpha)} \int_0^t (t-s)^{\alpha-1} N(s) \mathrm ds
					\right)^2 \right| N(t) = n \right\} \\
					= {} & \frac{1}{\Gamma^2(\alpha)} \int_0^t \int_0^s (t-s)^{\alpha-1} (t-w)^{\alpha-1}
					\mathbb{E}\left( \left. N(s)N(w) \right| N(t) = n \right) \mathrm ds \, \mathrm dw \notag \\
					& + \frac{1}{\Gamma^2(\alpha)} \int_0^t \int_s^t (t-s)^{\alpha-1} (t-w)^{\alpha-1}
					\mathbb{E}\left( \left. N(s)N(w) \right| N(t) = n \right) \mathrm ds \, \mathrm dw \notag \\
					= {} & \frac{n}{t \Gamma^2(\alpha)} \int_0^t (t-s)^{\alpha-1}\, \mathrm ds
					\int_0^s w(t-w)^{\alpha-1} \, \mathrm dw
					+	\frac{n}{t \Gamma^2(\alpha)} \int_0^t s(t-s)^{\alpha-1} \, \mathrm ds \int_s^t (t-w)^{\alpha-1}
					\, \mathrm dw \notag \\
					& + \frac{n(n-1)}{t^2 \Gamma^2(\alpha)} \int_0^t \int_0^t
					(t-s)^{\alpha-1} (t-w)^{\alpha-1} sw \, \mathrm ds \, \mathrm dw \notag \\
					= {} & \frac{2nt^{2\alpha} \Gamma(2\alpha)}{\alpha \Gamma^2(\alpha)\Gamma(2\alpha+2)}
					+ \frac{n(n-1)t^{2\alpha}}{\Gamma^2(\alpha+2)}. \notag
				\end{align}
				From results \eqref{03re} and \eqref{02blabla}, the conditional variance then simply reads
				\begin{align}
					\mathbb{V}\text{ar} \left\{ \left. \frac{1}{\Gamma(\alpha)} \int_0^t (t-s)^{\alpha-1} N(s)
					\, \mathrm ds
					\right| N(t) = n \right\}
					= \frac{2 n t^{2\alpha}\Gamma(2\alpha)}{\Gamma(\alpha)\Gamma(\alpha+1)\Gamma(2\alpha+2)}
					- \frac{nt^{2\alpha}}{\Gamma^2(\alpha+2)}
					= \frac{nt^{2\alpha} \alpha^2}{(2\alpha+1)\Gamma^2(\alpha+2)}.
				\end{align}
			\end{proof}
		\end{thm}
		
		\begin{rem}
			The unconditional variance \eqref{05cc} can be easily retrieved as follows.
			\begin{align}
				\label{06fin}
				& \mathbb{V}\text{ar} \left( \mathcal{N}^\alpha(t) \right) = \\
				& \mathbb{E} \left\{ \mathbb{V}\text{ar} \left( \mathcal{N}^\alpha(t) |N(t) \right) \right\}
				+ \mathbb{V}\text{ar} \left\{ \mathbb{E}\left( \mathcal{N}^\alpha(t) | N(t) \right) \right\}
				= \frac{\lambda t^{2\alpha+1}\alpha^2}{(2\alpha+1)\Gamma^2(\alpha+2)}
				+ \frac{\lambda t^{2\alpha+1}}{\Gamma^2(\alpha+2)}
				= \frac{\lambda t^{2\alpha+1}}{(2\alpha+1)\Gamma^2(\alpha+1)}. \notag
			\end{align}
		\end{rem}

	\section{Integral of the homogenous Poisson process}
	
		For the integral of the Poisson process we have a representation in terms of random sums.
		
		\begin{thm}
			For the homogenous Poisson process $N(t)$, $t \ge 0$, we have that
			\begin{align}
				\label{08icic}
				\mathcal{N}(t) = \int_0^t N(s) \, \mathrm ds \overset{\text{d}}{=}
				\sum_{j=1}^{N(t)} X_j = \mathfrak{N}(t),
			\end{align}
			where the $X_j$s are i.i.d.\ random variables, uniform in $[0,t]$. In \eqref{08icic}, the sum in the right
			hand side is intended to be equal to zero when $N(t)=0$.
			\begin{proof}
				If $N(t) = n$, and $\tau_1, \tau_2, \dots, \tau_n$, are the random instants at which the Poisson
				events appear, we have that the integral of the Poisson $A_n(t)$, $t \ge 0$, reads
				\begin{align}
					A_n(t) = \sum_{j=2}^n (\tau_j-\tau_{j-1})(j-1) + n(t-\tau_n).
				\end{align}
				Since
				\begin{align}
					\Pr \{ \tau_1 \in \mathrm d s_1, \dots, \tau_n \in \mathrm ds_n \} = \frac{n!}{t^n},
					\qquad 0 < s_1 < s_2 < \dots < s_n < t,
				\end{align}
				we have
				\begin{align}
					\mathbb{E} e^{i\beta A_n(t)} = \frac{n!}{t^n} \int_0^t \mathrm d{s_1} \dots
					\int_{s_{n-1}}^t \mathrm d{s_n} e^{i\beta \left[ \sum_{j=2}^n
					(s_j-s_{j-1}) (j-1) + n(t-s_n) \right]}
					= \frac{n!}{t^n} F_n(\beta,t).
				\end{align}
				It is evident that the functions $F_n(\beta, t)$ satisfy the equations
				\begin{align}
					\label{nestor}
					\begin{cases}
						\frac{\mathrm d}{\mathrm d{t}} F_n(\beta,t) = i \, n \, \beta F_n(\beta,t) + F_{n-1}(\beta,t),
						& n \geq 1, \\
						F_n(\beta,0) = 0,
					\end{cases}
				\end{align}
				where $F_0(\beta,t)=0$.
				Now we show by induction that
				$
					F_n(\beta,t) = \left( e^{i\beta t}-1 \right)^n / \left( n! (i\beta)^n \right)
				$, $n \geq 0$.
				From \eqref{nestor} we have that
				\begin{align}
					\frac{\mathrm d}{\mathrm d{t}} F_n(\beta,t) - i \, n \, \beta F_n(\beta,t) =
					\frac{\left( e^{i\beta t} -1 \right)^{n-1}}{\left(
					n-1 \right)! \left( i\beta \right)^{n-1}},
				\end{align}
				and the solutions turn out to be
				\begin{align}
					F_n(\beta,t) & = e^{in\beta t} \left[ \int_0^t \frac{\left( e^{i\beta s}-1
					\right)^{n-1}}{\left( n-1 \right)! \left( i\beta \right)^{n-1}}
					e^{-i n \beta s} \mathrm d{s} \right]
					= \frac{e^{in\beta t}}{\left( n-1 \right)! \left( i\beta \right)^{n-1}}
					\sum_{m=0}^{n-1} \binom{n-1}{m} (-1)^{n-1-m} \int_0^t e^{im\beta s - in\beta s}
					\mathrm d{s} \\
					& = \frac{e^{i n\beta t}}{(n-1)!(i\beta)^{n-1}} \frac{1}{i\beta}
					\sum_{m=0}^{n-1} \binom{n-1}{m} (-1)^{n-1-m} \frac{\left( e^{-i\beta t
					\left( n-m \right)} -1 \right)}{(m-n)} \notag \\
					& = \frac{e^{i n\beta t}}{(n-1)!(i\beta)^n} \sum_{m=0}^{n-1} \frac{(n-1)!}{
					m!(n-1-m)!(n-m)} (-1)^{n-m} \left[ e^{-i\beta t \left( n-m \right)} -1 \right] \notag \\
					& = \frac{e^{i\beta n t}}{n!(i\beta)^n} \left[ \sum_{r=0}^{n-1} (-1)^{n-r}
					\left( e^{-i\beta t\left( n-r \right)}-1 \right) \right]
					= \frac{e^{i\beta n t}}{n!(i\beta)^n} \left[ \sum_{r=1}^n \binom{n}{r}
					(-1)^r \left( e^{-i\beta t r}-1 \right) \right] \notag \\
					& = \frac{e^{i\beta n t}}{n!(i\beta)^n} \left[ \sum_{r=0}^n \binom{n}{r}
					(-1)^r \left( e^{-i\beta t r}-1 \right) \right]
					= \frac{e^{i\beta n t}\left( 1-e^{-i \beta t} \right)^n}{n! (i \beta)^n}
					= \frac{\left( e^{i\beta t} -1 \right)^n}{n! (i\beta)^n}. \notag
				\end{align}
				The characteristic function of $A_n(t)$ can thus be written as
				\begin{align}
					\mathbb{E} e^{i \beta A_n(t)} = \frac{\left( e^{i\beta t}
					-1 \right)^n}{t^n(i \beta)^n},
				\end{align}
				so that
				\begin{align}
					\label{05sa}
					\mathbb{E} e^{i \beta \mathcal{N}(t)}
					= e^{-\lambda t} \sum_{n=0}^\infty \left( e^{i\beta t} -1 \right)^n
					\frac{(\lambda t)^n}{t^n} \frac{1}{(i\beta)^n} \frac{1}{n!}
					= e^{-\lambda t + \frac{\lambda}{i\beta}\left( e^{i\beta t} -1 \right)}
					= e^{\lambda \int_0^t \left( e^{i\beta s} -1 \right) \mathrm d{s} },
				\end{align}
				which is the characteristic function of the compound process $\mathfrak{N}(t)$, $t\ge 0$.
			\end{proof}		
		\end{thm}
		
		\begin{rem}
			From \eqref{05sa}, we have that
			\begin{align}
				\mathbb{E} e^{i\beta \mathcal{N}(t)} = e^{i\lambda\frac{\beta t^2}{2} - \lambda \frac{\beta^2
				t^3}{6} + o(t^3)}.
			\end{align}
			This shows that for small $t$ the integrated Poisson process is Gaussian with mean
			$\lambda t^2/2$ and variance $\lambda t^3/3$. The parameters of the approximating Gaussian
			coincide with the mean and variance of $\mathcal{N}(t)$.
		\end{rem}
		
		In the previous section we have obtained that
		\begin{align}
			\mathbb{E} \left\{ \mathcal{N}^\alpha(t) | N(t) = n \right\} = \frac{nt^\alpha}{\Gamma(\alpha+2)},
			\qquad n \ge 0,
		\end{align}
		and thus, for $\alpha=1$,
		$
			\mathbb{E} \left\{ \mathcal{N}(t) | N(t) = n \right\}  = nt/2.
		$
		We are able to derive this result with a different technique and, in the same way, to obtain
		\begin{align}
			\mathbb{E} \left\{ \int_0^t \left[  N(s) \right]^k \mathrm ds \biggr| N(t) = n \right\}
		\end{align}
		for $k=2,3$. The same technique is applied for the derivation of
		\begin{align}
			\mathbb{V}\text{ar} \left\{ \int_0^t N(s) \, \mathrm ds \biggr| N(t) = n \right\}.
		\end{align}
		
		Before stating the next theorem we recall again that the conditional distribution of
		$N(s)$, given $N(t) = n$, $s < t$, is Binomial$(n,s/t)$ (see formula \eqref{labello} and
		\citet[page 21]{kingman}).
		\begin{thm}
			\label{06teo}
			For the integrated powers of the Poisson process we have that
			\begin{align}
				\mathbb{E} \left\{ \int_0^t \left[ N(s) \right]^k \mathrm ds \biggr| N(t) = n \right\}
				= \frac{t}{n+1} \sum_{j=1}^n j^k.
			\end{align}
			\begin{proof}
				For $0=t_0<t_1< \dots < t_n< t_{n+1}= t$, under the condition that $N(t) = n$, we can write
				\begin{align}
					\int_0^t \left[ N(s) \right]^k \mathrm ds = \sum_{j=1}^{n+1} \int_{t_{j-1}}^{t_j}
					\left[ N(s) \right]^k \mathrm ds
					= \sum_{j=1}^{n+1} (j-1)^k (t_j-t_{j-1}).
				\end{align}
				Therefore,
				\begin{align}
					& \mathbb{E} \left\{ \int_0^t \left[ N(s) \right]^k \mathrm ds
					\biggr| N(t)=n \right\}
					= \frac{n!}{t^n} \int_0^t \mathrm dt_1 \dots \int_{t_{j-1}}^t \mathrm dt_j
					\dots \int_{t_{n-1}}^t \mathrm dt_n \sum_{j=1}^{n+1} (j-1)^k (t_j-t_{j-1}) \\
					& = \frac{n!}{t^n} \sum_{j=1}^{n+1} \frac{(j-1)^k}{(n-j)!} \int_0^t \mathrm dt_1
					\dots \int_{t_{j-1}}^t (t_j-t_{j-1})(t-t_j)^{n-j} \mathrm dt_j. \notag
				\end{align}
				Since
				\begin{align}
					\int_{t_{j-1}}^t (t_j-t_{j-1})(t-t_j)^{n-j} \mathrm dt_j
					& = \int_0^{t-t_{j-1}} w^{n-j} (t-t_{j-1}-w)\, \mathrm dw
					& = (t-t_{j-1})^{n-j+2} \frac{\Gamma(2)\Gamma(n-j-1)}{\Gamma(n-j+3)},
				\end{align}
				we arrive at
				\begin{align}
					\mathbb{E} \left\{ \int_0^t \left[ N(s) \right]^k \mathrm ds
					\biggr| N(t)=n \right\}
					& = \frac{n!}{t^n} \sum_{j=1}^{n+1} \frac{(j-1)^k}{(n-j)!} \int_0^t \mathrm dt_1
					\dots \int_{t_{j-2}}^t \mathrm dt_{j-1} (t-t_{j-1})^{n-j+2} \frac{(n-j)!}{(n-j+2)!} \\
					& = \frac{n!}{t^n} \sum_{j=1}^{n+1} \frac{(j-1)^k}{(n+1)!}
					= \frac{t}{n+1} \sum_{j=1}^n j^k. \notag
				\end{align}
			\end{proof}
		\end{thm}
		
		\begin{rem}
			Explicit results can be given for small values of $k$:
			\begin{align}
				\mathbb{E} \left\{ \int_0^t \left[ N(s) \right]^k \mathrm ds \biggr| N(t)=n \right\} =
				\begin{cases}
					\frac{nt}{2}, & k=1, \\
					\frac{n(2n+1)t}{6}, & k=2, \\
					\frac{n^2(n+1)t}{4}, & k=3.
				\end{cases}
			\end{align}
			
			The unconditional mean values have therefore the form
			\begin{align}
				\mathbb{E} \left\{ \int_0^t \left[ N(s) \right]^k \mathrm ds \right\} =
				\begin{cases}
					\frac{\lambda t}{2}, & k=1, \\
					\frac{\lambda^2t^3}{3} + \frac{\lambda t^2}{2}, & k=2, \\
					\frac{\lambda^3 t^4}{4} + \lambda^2 t^3 + \frac{\lambda t^2}{2}, & k=3.
				\end{cases}
			\end{align}
		\end{rem}
		
		By applying the same technique as in Theorem \ref{06teo} we obtain the conditional variance.
		
		\begin{thm}
			We have the following explicit results.
			\begin{align}
				\label{06uno}
				\mathbb{E} \left\{ \left[ \int_0^t N(s)\, \mathrm ds \right]^2
				\biggr| N(t)= n \right\} = \frac{n(3n+1)t^2}{12}, \qquad
				\mathbb{V}\text{ar} \left\{ \int_0^t N(s) \, \mathrm ds \biggr| N(t) = n \right\}
				= \frac{nt^2}{12}.
			\end{align}
			\begin{proof}
				If we assume that $N(t) = n$, the following decomposition holds.
				\begin{align}
					\label{06cic}
					\left( \int_0^t N(s) \, \mathrm ds \right)^2
					= \sum_{j=1}^{n+1} (j-1)^2 (t_j-t_{j-1})^2
					+ 2 \sum_{1\le j< r \le n+1} (j-1)
					(r-1)(t_j-t_{j-1})(t_r-t_{r-1}),
				\end{align}
				for $0=t_0<t_1< \dots <t_j<\dots<t_n< t_{n+1}=t$. Note that
				\begin{align}
					\label{06aa}
					\frac{n!}{t^n} & \int_0^t \mathrm dt_1 \dots \int_{t_{j-1}}^t (t_j-t_{j-1}) \, \mathrm dt_j
					\int_{t_j}^t \mathrm dt_{j+1} \dots \int_{t_{r-1}}^t (t_r-t_{r-1})
					\, \mathrm dt_r
					\int_{t_r}^t \mathrm dt_{r+1} \dots \int_{t_{n-1}}^t \mathrm dt_n \\
					= {} & \frac{n!}{t^n} \int_0^t \mathrm dt_1 \dots \int_{t_{j-1}}^t (t_j-t_{j-1})
					\, \mathrm dt_j
					\int_{t_j}^t \mathrm dt_{j+1} \dots \int_{t_{r-1}}^t (t_r-t_{r-1})
					\frac{(t-t_r)^{n-r}}{(n-r)!} \mathrm dt_r \notag \\
					= {} & \frac{n!}{t^n \Gamma(n-r+3)} \int_0^t \mathrm dt_1 \dots
					\int_{t_{j-1}}^t (t_j-t_{j-1}) \, \mathrm dt_j
					\int_{t_j}^t \mathrm dt_{j+1}
					\int_{t_{r-2}}^t \mathrm dt_{r-1} (t-t_{r-1})^{n-r+2} \notag \\
					= {} & \frac{n!}{t^n \Gamma(n-j+2)} \int_0^t \mathrm dt_1 \dots \int_{t_{j-1}}^t
					\mathrm dt_j (t-t_j)^{n-j+1} (t_j-t_{j-1}) \notag \\
					= {} & \frac{n!}{t^n\Gamma(n-j+2)} \int_0^t \mathrm dt_1 \dots \int_{t_{j-2}}^t
					\mathrm dt_{j-1} (t-t_{j-1})^{n-j+3} \frac{\Gamma(n-j+2)\Gamma(2)}{\Gamma(n-j+4)} \notag \\
					= {} & \frac{n!}{\Gamma(n-j+4)} \frac{t^2}{(n-j+4)(n-j+5)\dots (n+2)} \notag
					= \frac{t^2}{(n+2)(n+1)}.
				\end{align}
				In the same way we have that
				\begin{align}
					\label{06beb}
					& \frac{n!}{t^n} \int_0^t \mathrm dt_1 \dots \int_{t_{j-1}}^t (t_j-t_{j-1})^2
					\mathrm d t_j \int_{t_j}^t \mathrm d t_{j+1} \dots \int_{t_{n-1}}^t \mathrm dt_n
					= \frac{n!}{t^n} \int_0^t \mathrm dt_1 \dots \int_{t_{j-1}}^t (t_j-t_{j-1})^2
					\frac{(t-t_j)^{n-j}}{(n-j)!} \mathrm dt_j \\
					& = \frac{n!}{t^n} \int_0^t \mathrm dt_1 \dots \int_{t_{j-2}}^t \mathrm dt_{j-1}
					(t-t_{j-1})^{n-j+3} \frac{2}{(n-j+3)!}
					= \frac{2n! t^{n+2}}{t^n (n+2)!}
					= \frac{2 t^2}{(n+1)(n+2)}. \notag
				\end{align}
				In light of formulae \eqref{06aa}, \eqref{06beb} and decomposition \eqref{06cic}, we have
				\begin{align}
					\label{nes2}
					\mathbb{E} & \left\{ \left( \int_0^t N(s) \mathrm ds \right)^2 \biggr| N(t) = n \right\}
					= \frac{2 t^2}{(n+1)(n+2)} \sum_{j=1}^{n+1} (j-1)^2 + 2 \sum_{j=1}^n (j-1)
					\sum_{r=j+1}^{n+1} (r-1) \frac{t^2}{(n+1)(n+2)} \\
					= {} & \frac{2 t^2}{(n+1)(n+2)} \left[ \sum_{j=1}^n j^2 + \sum_{j=1}^n (j-1)
					\sum_{r=j+1}^{n+1} (r-1) \right] \notag \\
					= {} & \frac{2 t^2}{(n+1)(n+2)} \left[ \frac{n(n+1)(2n+1)}{6} + \sum_{j=2}^n
					(j-1) \left( \frac{n(n+1)}{2} - \frac{j(j-1)}{2} \right) \right] \notag \\
					= {} & \frac{2 t^2}{(n+1)(n+2)} \left[ \frac{n(n+1)(2n+1)}{6} + \frac{n^2(n+1)(n-1)}{4} \right.
					\left. -\frac{1}{2} \left( \frac{n(n-1)}{2} \right)^2 -\frac{1}{2}
					\frac{n(n-1)(2(n-1)+1)}{6} \right] \notag \\
					= {} & \frac{t^2}{(n+1)(n+2)} \left[ \frac{n(n+1)(2n+1)}{3} - \frac{n(n-1)(2n-1)}{6} \right.
					\left. + \frac{n^2(n-1)(n+3)}{4} \right] \notag \\
					= {} & \frac{n t^2}{12 (n+1)(n+2)} \left[ 4(n+1)(2n+1) \right.
					\left. - 2(n-1)(2n-1) + 3n(n-1)(n+3) \right] \notag \\
					= {} & \frac{nt^2}{12(n+2)} \left[ 4(2n+1) + (n-1)(3n+2) \right]
					= \frac{n(3n+1)t^2}{12}. \notag
				\end{align}
				The conditional variance easily follows.
			\end{proof}
		\end{thm}
		
		\begin{rem}
			The results of Theorem \ref{06teorema}, for $\alpha=1$, coincide with \eqref{06uno}.
			
			We also observe that
			\begin{align}
				\mathbb{V}\text{ar} \left\{ \int_0^t N(s) \, \mathrm ds \right\} =
				\mathbb{E} \left\{ \mathbb{V}\text{ar} \left( \int_0^t N(s) \, \mathrm ds \biggr| N(t) \right) \right\}
				+ \mathbb{V}\text{ar} \left\{ \mathbb{E} \left( \int_0^t N(s) \, \mathrm ds
				\biggr| N(t) \right) \right\}
				= \frac{\lambda t^3}{3},
			\end{align}
			and this coincides with \eqref{06fin} for $\alpha=1$.
		\end{rem}
		
		\subsection{Properties of the integral of $\mathring{N}(t) = N_\lambda(t)-N_\beta(t)$}
		
			\label{prop99}
			It is well-known that for two independent Poisson processes the difference
			$\mathring{N}(t) = N_\lambda(t)-N_\beta(t)$, $t \ge 0$, (which can be used, for example, in modelling
			immigration-emigration processes) has Skellam distribution:
			\begin{align}
				\Pr \{ \mathring{N}(t) = r  \} = e^{-(\beta+\lambda)t} \left( \frac{\lambda}{\beta} \right)^{r/2}
				I_{|r|}(2 t \sqrt{\lambda \beta}), \qquad r \in \mathbb{Z}, \: t>0,
			\end{align}
			where
			\begin{align}
				I_\alpha(z) = \sum_{m=0}^\infty \frac{(z/2)^{2m+\alpha}}{m!\Gamma(m+\alpha+1)}
			\end{align}
			is the modified Bessel function of the first kind.
			
			For the integral of the difference of the two Poisson processes we have that
			\begin{align}
				\label{baba}
				& \mathbb{E} e^{i\mu \left( \int_0^t N_\lambda(s) \mathrm ds
				- \int_0^t N_\beta(s) \mathrm ds \right)} = e^{\lambda \int_0^t \left( e^{i\mu s} -1 \right)
				\mathrm ds + \beta \int_0^t \left( e^{-i\mu s} -1 \right) \mathrm ds}
				= e^{-(\lambda+\beta)t +\int_0^t \left( \lambda e^{i\mu s} + \beta e^{-i\mu s} \right) \mathrm ds}
				\\
				& = e^{-(\lambda+\beta)t + \frac{\lambda e^{i\mu t}-\lambda}{i\mu} - \frac{\beta e^{-i\mu t}-\beta}{i\mu}}
				= e^{-(\lambda+\beta)t +\frac{1}{i\mu}(\lambda-\beta)\cos \mu t + \frac{(\lambda+\beta)}{\mu}
				\sin \mu t - \frac{(\lambda -\beta)}{i\mu}}
				= e^{(\lambda + \beta)t \left( \frac{\sin \mu t}{\mu t} -1 \right)}
				e^{-\frac{(\lambda-\beta)}{i\mu}\left( 1-\cos \mu t \right)}. \notag
			\end{align}

			\begin{thm}
				For the difference of integrated Poisson processes we have that
				\begin{align}
					\int_0^t \mathring{N}(s)\mathrm ds =
					\int_0^t N_\lambda(s) \mathrm ds - \int_0^t N_\beta (s) \mathrm ds \overset{\mathrm d}{=}
					\sum_{j=1}^{\widetilde{N}(t)} Z_j,
				\end{align}
				where $\widetilde{N}(t)$, $t \ge 0$,
				is a Poisson process of rate $\lambda+\beta$ and the $Z_j$s are i.i.d.\
				random variables with law
				\begin{align}
					f(s) =
					\begin{cases}
						\frac{\beta}{t(\lambda+\beta)}, & -t < s \le 0, \\
						\frac{\lambda}{t(\lambda+\beta)}, & 0 < s < t.
					\end{cases}
				\end{align}
				\begin{proof}
					The claimed result can be proved by resorting to the characteristic function.
					\begin{align}
						\label{07unalabel}
						\mathbb{E} e^{i \mu \left[ \int_0^t N_\lambda(s) \mathrm ds - \int_0^t N_\beta (s)
						\mathrm ds \right]}
						= e^{\lambda \int_0^t (e^{i \mu s}-1) \mathrm ds + \beta \int_0^t (e^{-i\mu s}-1) \mathrm ds}
						= e^{(\lambda+\beta) t \int_{-t}^t \left(e^{i\mu s}-1 \right)
						\left( \frac{\lambda}{(\lambda+\beta)t}
						\mathbb{I}_{[0,t]}(s) + \frac{\beta}{(\lambda + \beta)t} \mathbb{I}_{[-t,0]}(s) \right) 
						\mathrm ds}.
					\end{align}
				\end{proof}
			\end{thm}
			
			\begin{rem}
				From \eqref{07unalabel} we see that, for small values of $t$, the integrated difference of Poisson
				processes has Gaussian distribution with mean $(\lambda-\beta)t^2/2$ and variance $(\lambda+\beta)t^3/3$.
			\end{rem}
			
	\bibliographystyle{plainnat}
	\bibliography{bjps}
	\nocite{*}

\end{document}